\documentclass[openacc]{rstransa}
\pdfoutput=1 
\usepackage{float}
\usepackage{subfig}

\usepackage[linesnumbered,ruled,vlined]{algorithm2e}
\usepackage{siunitx}
\SetKwInput{KwInit}{Input}
\usepackage{mathtools}

\newcommand{\ve}[1]{\pmb{#1}}

\newcommand{\bu}{\ve{u}}

\newcommand{\bw}{\ve{w}}

\newcommand{\bI}{\ve{I}}

\newcommand{\bxi}{\ve{\xi}}

\newcommand{\bnabla}{\ve{\nabla}}

\newcommand{\mom}{\mu}
\newcommand{\bmom}{\ve{\mu}}

\newcommand\Rey{\mathrm{Re}}

\newcommand{\Rdot}{\dot{R}}
\newcommand{\Rddot}{\ddot{R}}

\begin{document}

\title{Hybrid quadrature moment method for accurate and stable representation of non-Gaussian processes and their dynamics}

\author{A. Charalampopoulos$^{1}$, S. H. Bryngelson$^{2}$, T. Colonius$^{3}$, and T. P. Sapsis$^{1}$}

\address{$^{1}$Department of Mechanical Engineering, Massachusetts Institute of Technology, Cambridge, MA 02139, USA\\
$^{2}$School of Computational Science and Engineering, Georgia Institute of Technology, GA 30313, USA\\
$^{3}$Division of Engineering and Applied Science, California Institute of Technology, Pasadena, CA 91125, USA}

\subject{
computational physics, computer modelling and simulation, mechanics
}

\corres{A. Charalampopoulos\\
\email{alexchar@mit.edu}
}

\begin{abstract}
Solving the population balance equation (PBE) for the dynamics of a dispersed phase coupled to a continuous fluid is expensive. Still, one can reduce the cost by representing the evolving particle density function in terms of its moments.  
In particular, quadrature-based moment methods (QBMMs) invert these moments with a quadrature rule, approximating the required statistics.
QBMMs have been shown to accurately model sprays and soot with a relatively compact set of moments. However, significantly non-Gaussian processes such as bubble dynamics lead to numerical instabilities when extending their moment sets accordingly.
We solve this problem by training a recurrent neural network (RNN) that adjusts the QBMM quadrature to evaluate unclosed moments with higher accuracy.
The proposed method is tested on a simple model of bubbles oscillating in response to a temporally fluctuating pressure field.
The approach decreases model-form error by a factor of $10$ when compared to traditional QBMMs. It is both numerically stable and computationally efficient since it does not expand the baseline moment set.
Additional quadrature points are also assessed, optimally placed and weighted according to an additional RNN. 
These points further decrease the error at low cost since the moment set is again unchanged.
\end{abstract}


\begin{fmtext}


\end{fmtext}


\maketitle

\section{Introduction}\label{s:intro}

The dynamics of dispersions of small particles or bubbles in a fluid are important to many engineering and medical applications. 
In medicine, ultrasounds, generated via small cavitating bubbles, are employed during cataract removal~\cite{kelman1967phaco}, to stop internal bleeding~\cite{vaezy1997liver, vaezy1999hemostasis}, and in other procedures like tumor necrosis~\cite{bailey2003physical}.
Focused shockwaves can cavitate bubbles that ablate kidney stones during lithotripsy treatment~\cite{coleman1987acoustic, pishchalnikov2003cavitation}. 
Their interaction with biological tissue or manufactured soft materials also attracts the medical~\cite{brennen2015cavitation, pan2018phenomenology, oguri2018cavitation, dollet2019bubble, spratt2021characterizing} and material science communities~\cite{gaudron2015bubble, barajas2017effects, turangan2017numerical}.
Bubble cavitation is also responsible for damage and noise in hydraulic pipe systems~\cite{weyler1971investigation, streeter1983transient}, hydro turbines~\cite{escaler2006detection, kumar2010study, luo2016review}, and propellers~\cite{sharma1990cavitation, ji2012numerical}.
At the same time, soots are critical to combustion~\cite{kazakov1998dynamic, balthasar2003stochastic, pedel2014large, mueller2009joint} and aerosols are used in many industrial processes~\cite{sibra2017simulation, laurent2001multi, hussain2015new}. 
In nature, cavitation is used as part of the hunting strategies of some marine animals, including humpback whales~\cite{leighton2004trapped, leighton2007acoustical, bryngelson2020simulation}, mantis shrimps~\cite{patek2004deadly}, and snapping shrimps~\cite{bauer2004remarkable, koukouvinis2017unveiling}.

While the dynamics of these particles can be simulated directly for a specific (sampled) dispersion by tracking each particle, distribution statistics are typically sought in applications.  
In flows with large spatial gradients, a large ensemble of such simulations (Monte Carlo) is required to gather these statistics~\cite{zhao2007analysis, rosner2003multivariate}.  
The poor scaling of Monte Carlo makes such simulations expensive, and particle tracking also interferes with efficient parallelization.  
By instead phase-averaging the equations of motion~\cite{zhang1994ensemble}, a two-way coupled set of Eulerian equations that are more suitable to parallelization and GPU processing is obtained.  
However, the averaged equations involve solving the generalized population balance equation (PBE)~\cite{ramkrishna2000population}. 
The PBE evolves the dispersed phase number density function (NDF) as a function of its dynamic variables~\cite{bryngelson2020qbmmlib}.
For example, the relevant variables for bubbles dynamics are the bubble radii and their radial velocities.
However, further treatment is still required.
The PBE is a partial differential equation in the dynamic variables, separate from the spatial and temporal variables of the flow equations, making this approach intractable for large simulations.

Quadrature-based moment methods (QBMMs) are a low-cost approach to approximately solving a PBE.
Introduced in~\cite{mcgraw1997description}, QBMMs have seen rapid improvement~\cite{marchisio2013computational}. 
In brief, they prescribe a finite moment set and invert it to an optimal set of quadrature nodes and weights in the dynamic system phase space. 
The success of QBMMs has led to the creation of open-source libraries for them~\cite{bryngelson2020qbmmlib, passalacqua2018open}.
In the case of multiple dynamic variables, conditional QBMMs like conditional-QMOM (CQMOM)~\cite{yuan2011conditional} and conditional hyperbolic-MOM (CHyQMOM)~\cite{fox2018conditional, patel2019three} are preferred.
These methods can efficiently solve many problems but suffer from a combinatorial explosion of their computational cost when higher accuracy is needed. 
This problem stems from the need to evolve all moments up to a higher order to increase accuracy. 
Worse still, these methods can exhibit numerical instabilities when third- or higher-order moments are evolved~\cite{marchisio2013computational}.

To circumvent these stability issues, this work employs neural networks to enhance the predictive abilities of standard $2$-by-$2$-node ($4$-node) CHyQMOM, which only requires access to first and second-order moments. 
This approach avoids both the numerical instabilities and high computational costs of evolving higher-order moments.
The method follows the recent success of deep neural networks for improving multiphase flow models~\cite{ma2015using, ma2016using, charalampopoulos2021machine, bryngelson2020gaussian}.
We expand on a previous such effort that used neural networks to close strictly-Gaussian moment transport equations~\cite{bryngelson2020gaussian}.
Here, we instead seek data-informed corrections to a CHyQMOM method~\cite{fox2018conditional, patel2019three}.
By doing this, one has control over the resulting quadrature nodes and weights.
This makes correcting moment approximations straightforward and consolidates the two neural networks of~\cite{bryngelson2020gaussian} to one.
This allows for computation of even out-of-training-set moments, in contrast to data-informed moment methods that use low-order moments to learn specific high-order ones~\cite{huang2021machinea, huang2021machineb, huang2021machinec}.
Extension from~\cite{bryngelson2020gaussian} includes non-uniform and long-time pressure forcings, making the trained model appropriate for CFD solvers. 

We emphasize that incorporating a trained neural network into a numerical method comes with its constraints. 
The neural networks train on a finite set of pressure profiles thought appropriate for a class of physical problems. However, one can not ensure that this fully generalizes, something which must be verified.
Here, pressure profiles are sampled from a distribution that can represent many practical bubbly flow problems. 
Yet, the method is not guaranteed to generalize well for drastically different external forcings.

Section~\ref{s:formulation} formulates a model problem that serves as the basis for the proposed extension of CHyQMOM, namely the dynamics of a population of cavitating bubbles whose statistics are significantly non-Gaussian. 
In section~\ref{s:hybrid}, the new hybrid method is described. 
Section~\ref{s:results} shows that this approach can improve low-order moment predictions while extrapolating out of the moment space to compute required high-order moment predictions. 
This section also investigates the utility of additional quadrature points whose locations are selected by the RNN.  
It also compares the computational costs of the present approach and the classical CHyQMOM.
Section~\ref{s:conclusions} summarizes our conclusions.
      
\section{Problem Formulation}\label{s:formulation}

\subsection{Ensemble-averaged flow equations}

This work focuses on the fluid-coupled dynamics of a dispersion of small, spherical bubbles transported in a compressible carrier fluid. 
The mixture phase-averaged evolution equations for the continuous fluid are
\begin{gather}\label{Mixture_Eqs}
\begin{split}
    \frac{\partial \rho}{\partial t} + \bnabla \cdot (\rho \bu) &= 0, \\
    \frac{\partial \rho \bu}{\partial t} +\bnabla \cdot (\rho \bu \bu +p \bI) &= 0,\\
    \frac{\partial E}{\partial t} +\bnabla \cdot (E +p) \bu &= 0,
\end{split}
\end{gather}
where $\rho, \bu, p$, and $E$ being the mixture density, velocity vector, pressure, and total energy, respectively. 
The system of equations is complemented by appropriate initial and boundary or radiation conditions specific to each individual problem. 
The void fraction of the bubbles is $\alpha$ and a dilute assumption $\alpha \ll 1$ is made.
The bubbles are defined by their instantaneous bubble radii $R$, its time derivative $\Rdot$.
The bubbles are assumed monodisperse and so have the same equilibrium radius $R_o$.

The mixture pressure $p$ is deduced from the ensemble phase-averaging method~\cite{zhang1994ensemble, bryngelson2019quantitative} as 
\begin{gather}\label{Mixture_Pres}
    p = (1-\alpha)p_l +\alpha \left( \frac{ \overline{R^3 p_{bw}} }{ \overline{R^3} } -\rho \frac{ \overline{R^3 \Rdot^2}   }{ \overline{R^3}  }  \right),
\end{gather}
where $p_{bw}$ and $p_l$ are the bubble wall and liquid pressures, respectively~\cite{ando2011numerical}. 
Liquid pressure $p_l$ follows from the stiffened-gas equation of state~\cite{menikoff1989riemann}, though this model can be substituted for another if required. 
The usual coefficients for water are used~\cite{bryngelson2019quantitative}.

The overbars in~\ref{Mixture_Pres} denote raw moments $\mom$ of the bubble dispersion as
\begin{gather}\label{Moments_Def}
    \overline{R^i \Rdot^j} = \mom_{i,j} = \int_{\Omega} R^i \Rdot^j f(R,\Rdot;R_o) \mathrm{d}R \mathrm{d}\Rdot,
\end{gather}
where $f$ is the number density function of the bubbles.
This paper focused on a new, improved method for computing these moments, which will be introduced in section~\ref{s:hybrid}.

The void fraction transports as~\cite{bryngelson2019quantitative}:
\begin{gather}\label{void_frac_Eq}
    \frac{\partial \alpha}{\partial t} +\bu \cdot \bnabla \alpha = 3 \alpha \frac{\overline{R^2 \Rdot}}{ \overline{R^3} }.
\end{gather}
The moments required to close the governing flow equations are thus
\begin{gather}\label{Target_Moments}
    \bmom_{\mathrm{Targ.}} = \left\{ \overline{R^3 \Rdot^2}, \overline{R^3}, \overline{R^2 \Rdot}, \overline{R^3 p_{bw}} \right\}.
\end{gather}

\subsection{Bubble model}
To close the governing equations of the previous subsection, a model for the bubble dynamics, in terms of the dynamical variables $R$ and $\Rdot$, is required.
We use a Rayleigh--Plesset equation for this:
\begin{gather}\label{RP_Equation}
\begin{split}
    R \Rddot + \frac{3}{2} \Rdot^2 +\frac{4}{\Rey} \frac{\Rdot}{R} = 
    \left( \frac{R_o}{R} \right)^{3 \gamma} -
    \frac{1}{C_p},
\end{split}
\end{gather}
which is dimensionless via the reference bubble size $R_o$, liquid pressure $p_0$, and liquid density $\rho_0$. 
In~\eqref{RP_Equation}, $C_p$ is the ratio between the fluid and bubble pressures and $\Rey$ is a Reynolds number
\begin{gather}
\begin{split}
    \Rey = \sqrt{\frac{p_0}{\rho_0}} \frac{R_o}{\nu_0}, 
\end{split}
\end{gather}
where $\nu_0$ is the liquid kinematic viscosity. 
For the cases considered here $\Rey = 10^3$.

This model assumes the bubbles remain spherical and compress via a polytropic process with coefficient $\gamma = 1.4$.
While this model can be generalized to include heat exchange and liquid compressibility, these effects are not critical to our study and thus omitted here.
Based on this model, the bubble wall pressure $p_{bw}$ simplifies the last moment of $\bmom_\mathrm{Targ.}$ as 
\begin{gather}
    \overline{R^3 p_{bw}} = \mom_{3(1-\gamma),0}.
\end{gather}

We also define a dimensionless time $t^\ast = t \omega_0$, where $\omega_0$ is the natural frequency of the bubbles.
To simplify the notation, $t$ will be used in place of $t^\ast$ hereon.


\subsection{Population balance formulation}

A number density function $f$  describes the statistics of the bubbles.
The generalized population balance equation is
\begin{gather}\label{PBE_Equation}
\begin{split}
    \frac{\partial f}{\partial t} +\frac{\partial }{ \partial R} ( f \Rdot) +\frac{\partial}{\partial \Rdot} (f \Rddot) = 0,
\end{split}
\end{gather}
assuming the bubbles do not coalesce or break up, though these effects can be included via empirically modeled terms if desired. 
QBMMs solve~\eqref{PBE_Equation} by representing $f$ as a set of raw moments $\bmom$~\cite{fox2018conditional, bryngelson2020qbmmlib}. 
Through an appropriate inversion procedure, these methods can transform these moments into quadrature nodes and weights in phase space. 
This allows for the approximation of $f$ via a weighted sum of Dirac delta functions. 
Hence, the following quadrature rule can approximate any raw moment 
\begin{gather}\label{Moments_QMOM}
    \mom_{i,j} = \sum_k w_{k} \xi^{i}_{1,k} \xi_{2,k}^j,
\end{gather}
where $\xi_{i,j}$ is the $j$-th quadrature point locations for the $i$-th internal coordinate.

\section{Hybrid quadrature moment method formulation}\label{s:hybrid}

We now present the hybrid, data-informed method for predicting the moments of cavitating bubbles.
Subsection~\ref{ss:chyqmom} presents $4$-node CHyQMOM predictions~\cite{bryngelson2020qbmmlib}.
Subsection~\ref{ss:correction} details the hybrid neural-network model that improves the predictions.

\subsection{CHyQMOM}\label{ss:chyqmom}

For two-dynamic-variable cases, conditioned moment methods are computationally preferable to traditional QMOM~\cite{yuan2011conditional}.
We use CHyQMOM because it can close a second-order moment system with fewer carried moments than CQMOM~\cite{fox2018conditional}.
For a $4$-node quadrature rule, it uses the first- and second-order moments
\begin{gather}\label{low_moments}
    \bmom = \{ \mom_{1,0}, \mom_{0,1}, \mom_{2,0}, \mom_{1,1}, \mom_{0,2} \}.
\end{gather}
The CHyQMOM inversion process for obtaining the nodes $\bxi$ and weights $\bw$ is presented in appendix~\ref{appendix:inversion}. 

Taking the time-derivative of each of the~\eqref{low_moments} moments and applying~\eqref{RP_Equation} results in
\begin{gather}\label{Moment_Equations}
\begin{split}
    \frac{\partial \mom_{1,0} }{\partial t} &= \mom_{0,1},\\
    \frac{\partial \mom_{0,1} }{\partial t} &= -\frac{3}{2} \mom_{-1,2} -\frac{4}{\Rey} \mom_{-2,1} +\mom_{-4,0} -C_p \mom_{-1,0},\\
    \frac{\partial \mom_{2,0}}{\partial t} &= 2 \mom_{1,1},\\
    \frac{\partial \mom_{1,1}}{\partial t} &= -\frac{1}{2} \mom_{0,2} -\frac{4}{\Rey} \mom_{-1,1} +\mom_{-3,0} -C_p,\\
    \frac{\partial \mom_{0,2}}{\partial t} &= -3 \mom_{-1,3} -\frac{8}{\Rey} \mom_{-2,2} +2 \mom_{-4,1} -2 C_p \mom_{-1,1},
\end{split}
\end{gather}
which are called the moment transport equations.
The quadrature rule~\eqref{Moments_QMOM} approximates unclosed moments in~\eqref{Moment_Equations}.

While this scheme is computationally cheap, it is challenging to extend to include additional quadrature points without high computational cost or numerical instabilities.
Thus, truncation errors can affect approximation of the right-hand-side of~\eqref{Moment_Equations} and the extrapolation out of the low-order moment space to the moments of $\bmom_{\mathrm{Targ.}}$ of~\eqref{Target_Moments}. 
We next present an augmented method that treats these issues without introducing numerical instability or high computational cost.

\subsection{Data-informed corrections}\label{ss:correction}

We improve the CHyQMOM moment inversion procedure by adding a correction term to the $4$-node quadrature rule and introducing additional quadrature nodes.
The unaugmented CHyQMOM quadrature rule is denoted via $\{ \bw^{\mathrm{(QBMM)}}, \bxi^{\mathrm{(QBMM)}} \}$.

For these corrections, an LSTM RNN is employed.
The LSTM incorporates non-Markovian memory effects into the reduced-order model. 
This approach is known to be capable of improving predictions of reduced-order models~\cite{vlachas2018data, bryngelson2020gaussian}. 

The corrections $\{ \bw^\prime, \bxi^\prime \}$ serve as input predictions for the first- and second-order moments as well as the pressure $\{ \mom_{1,0}, \mom_{0,1}, \mom_{2,0}, \mom_{1,1}, \mom_{0,2}, C_p \}$. 
They are modelled as
\begin{gather}
    \{ \bw^\prime(t), \bxi^\prime(t)  \}  = 
    \mathbb{G}[ \Theta; \bmom(\chi(t)), C_p(\chi(t)), \Rey],
\end{gather}
where the vector $\Theta$ denotes hyperparameters and optimized parameters of the neural network as obtained during training. 
More detail on the implementation is in section~\ref{s:results}, subsection~\ref{subs:tes_conf}. The chosen values for the hyperparameters are included in appendix~\ref{appendix:hyperparams}.
The history of the reduced-order model states is
\begin{gather}
    \chi(t) = \{ t, t -\tau_1, ..., t-\tau_N  \}.
\end{gather}
The hybrid quadrature rule follows as
\begin{gather}
    \bw = \bw^{\mathrm{(QBMM)}} + \bw^\prime 
    \quad \text{and}
    \quad \bxi = \bxi^{\mathrm{(QBMM)}}+\bxi^\prime. 
\end{gather}

The neural network loss function is designed to ensure the target high-order moments $\bmom_T$ can be accurately computed and that the low-order moments $\bmom$ evolve accurately.
Hence, the right-hand-side of~\eqref{Moment_Equations} is included in the loss function as
\begin{gather}\label{Loss_Function}
\begin{split}
    \mathcal{L} = \sum_{0 \leq i,j \leq 2} \alpha_{i,j} \left( \frac{\partial \mom_{i,j}^{\mathrm{(ML)}} }{\partial t} -\frac{\partial \mom_{i,j}^{\mathrm{(MC)}}}{\partial t} \right)^2+\sum_{(i,j) \in \mathcal{I}} \beta_{i,j} \left( \sum_k w_k \xi_{1,k}^i \xi_{2,k}^j -\mom_{i,j}^{\mathrm{(MC)}} \right)^2 \\
    +\lambda \sum_k \mathrm{Relu}(-w_k),
\end{split}
\end{gather}
where,
\begin{gather}
    \mathcal{I} = \{ (0,0), (1,0), (0,1), (2,0), (1,1), (0,2), (3,0), (2,1), (3,2), (3-3\gamma,0)  \},
\end{gather}
and 
\begin{gather}
    \alpha_{i,j} = \Bigg\lVert \frac{\partial \mom_{i,j}^{\mathrm{(MC)}}}{\partial t} \Bigg\rVert_{\infty}^{-1}, \quad 
    \beta_{i,j} = \Big\lVert \mom_{i,j}^{\mathrm{(MC)}} \Big\rVert_{\infty}^{-1}.
\end{gather}
The first term in~\eqref{Loss_Function} minimizes, in the $L^2$ sense, the right-hand-sides of~\eqref{Moment_Equations} (given $\bmom$).
The second term in~\eqref{Loss_Function} minimizes prediction error for both $\bmom$ and $\bmom_{\mathrm{Targ.}}$, while it also penalizes the network when the weights don't sum up to $1$ (under the assumption that $\mu_{0,0} = 1$). 
The last term in~\eqref{Loss_Function} penalizes negative-valued weights.

The discretized moment transport equation~\eqref{Moment_Equations} and the quadrature rule~\ref{Moments_QMOM} compute the time-derivatives $\partial \mom_{i,j}^{\mathrm{(ML)}}/ \partial t$ required in~\eqref{Loss_Function} as
\begin{gather}
\begin{split}
    \frac{\partial \mom_{1,0}^{\mathrm{(ML)}} }{\partial t} &= \sum_k w_k \xi_{1,k}\\
    \frac{\partial \mom_{0,1}^{\mathrm{(ML)}} }{\partial t} &= -\frac{3}{2} \sum_{k} w_k \xi_{1,k}^{-1} \xi_{2,k}^2 -\frac{4}{\Rey} \sum_{k} w_k \xi_{1,k}^{-2} \xi_{2,k} +\sum_{k} w_k \xi_{1,k}^{-4} -C_p \sum_{k} w_k \xi_{1,k}^{-1},\\
    \frac{\partial \mom_{2,0}^{\mathrm{(ML)}}}{\partial t} &= 2 \sum_{k} w_k \xi_{1,k} \xi_{2,k}\\
    \frac{\partial \mom_{1,1}^{\mathrm{(ML)}}}{\partial t} &= -\frac{1}{2} \sum_{k} w_k \xi_{2,k}^2 -\frac{4}{\Rey} \sum_{k} w_k \xi_{1,k}^{-1} \xi_{2,k} +\sum_{k} w_k \xi_{1,k}^{-3} - C_p,\\
    \frac{\partial \mom_{0,2}^{\mathrm{(ML)}}}{\partial t} &= -3 \sum_{k} w_k \xi_{1,k}^{-1} \xi_{2,k}^3 -\frac{8}{\Rey} \sum_{k} w_k \xi_{1,k}^{-2} \xi_{2,k}^2 +2 \sum_{k} w_k \xi_{1,k}^{-4} \xi_{2,k} -2 C_p \sum_{k} w_k \xi_{1,k}^{-1} \xi_{2,k}.
\end{split}
\end{gather}


Once trained, the scheme
\begin{gather*}
    \bmom \xrightarrow{\text{(ML)}} \{ \bw, \bxi \}    
\end{gather*}
results in a new quadrature rule that evaluates the right-hand-side of~\eqref{Moment_Equations}.
The moment transport equations~\eqref{Moment_Equations} then evolve via an adaptive $4$th-order Runge--Kutta time-stepper.
Algorithm~\ref{alg:RK4} describes this procedure.

\begin{algorithm}[H]\label{alg:RK4}
\SetAlgoLined
\KwInit{$\bmom = \{ \mom_{1,0}, \mom_{0,1}, \mom_{2,0}, \mom_{1,1}, \mom_{0,2}\}$;   $C_p, \Rey$.} 
\KwData{NN architecture, CHyQMOM method, 4th-order-accurate Runge--Kutta (RK4), 
$\Rey$, $C_p$, time interval $t \in [0,T]$, error-tolerance $\tau_\mathrm{tol}$, maximum time-step $\delta t_{\mathrm{max}}$.}
\KwResult{$\bmom(t_i)$ and $\bmom_\mathrm{Targ.}(t_i)$ for $i = 0, 1, \dots, n$}
Train $\mathbb{D}_{\bu}^{\mathrm{(ML)}}$ with $\mathbb{D}_{\bu}^{\mathrm{(MC)}}$\;
$n \leftarrow 0$ \;
  \While{$t \leq T$}{ 
  $s \leftarrow t$\;
  $ \Big\{ \bw^{\mathrm{(QBMM)}}, \bxi^{\mathrm{(QBMM)}}  \Big\} = \mathrm{CHyQMOM}\Big[\bmom(s)\Big]$ \tcp*{Moment inversion}
  $\{\bw^\prime, \bxi^\prime\}(s) = \mathbb{G} \Big[ \bmom(s); \Big\{ \bw^{\mathrm{(QBMM)}}, \bxi^{\mathrm{(QBMM)}}  \Big\}(s), C_p(s)  \Big]  $ \tcp*{ML correction}
  $ \{ \bw, \bxi  \}(s) = \Big\{ \bw^{\mathrm{(QBMM)}} +\bw^\prime, \bxi^{\mathrm{(QBMM)}} +\bxi^\prime  \Big\}(s)  $ \tcp*{Set quadrature rule}
  $\{ \bmom, \bmom_T, \partial \bmom/\partial t  \}(s) = \mathrm{Quadrature}\Big[\{ \bw, \bxi \}(s) \Big] $ \tcp*{Project to moment space}
  $\delta t \leftarrow \delta t_{\max}$\;
  $\mathrm{flag} \leftarrow 1$\;
  \While{$\mathrm{flag} > 1$}{
  $\delta t \leftarrow \delta t/2$\; 
  $\bmom^1(s+\delta t_{\max}) = \mathrm{RK4} \Big[ \{ \bmom, \partial \bmom/\partial t  \}(s); \delta t  \Big]$ \tcp*{Evolve low-order moments}
  $\bmom^2(s+\delta t_{\max}) = \mathrm{RK4} \Big[ \{ \bmom, \partial \bmom/\partial t  \}(s); \delta t/2  \Big] $ \tcp*{Evolve low-order moments}
  $\mathrm{flag} \leftarrow \mathrm{Floor}\left[ \max\limits_{0 \leq l+m \leq 2}{ \lVert \mu_{l,m}^1-\mu_{l,m}^2 \rVert_2/\tau_\mathrm{tol} } \right]$
  }
  $t \leftarrow s +\delta t$\;
  $n = n + 1$ \;
  }
 \caption{Hybrid CHyQMOM}
\end{algorithm}

Note that the closure terms need to be evaluated at times $t$, $t +\delta t/2$, and $t +\delta t$.
The neural network does not make predictions at $t+\delta/2$, so the equations are instead integrated in time by $2 \delta t$ instead of $\delta t$.

\section{Results}\label{s:results}

\subsection{Pressure signals}\label{subs:tes_conf}

The capabilities of the data-enhanced CHyQMOM method to predict the statistics of bubble populations are explored. 
The pressure term $C_p$ excites the bubbles causing oscillations. 
The representation of $C_p$ used here should be general enough to include pressure profiles seen in actual fluid flows. 
In a generic framework, let $C_p$ have a finite Fourier-series expansion
\begin{gather}\label{CP_Forcing}
\begin{split}
    C_p(t) = 1+\sum_{i=1}^{N} \alpha_i \sin \left[2 \pi f_i t+\phi_i \right],
\end{split}
\end{gather}
where $t$ corresponds to nondimensional time, nondimenionalized by the natural oscillation frequency of the bubbles $\omega_0$, $f_i$ are the included dimensionless frequencies, and $\alpha_i$ and $\phi_i$ are the corresponding amplitude and phase. 
It is stressed that $C_p = 1$ corresponds to the equilibrium pressure of the bubbles (for which $R=1$ and $\dot{R} = 0$).

Most cavitating flow applications do not contain pressure frequencies higher than the natural oscillation frequency of the bubbles~\cite{brennen2014cavitation} (with a notable exception of some high-intensity focused ultrasound treatments).
We operate under this constraint hereon, though higher frequencies could be included if desired.
On the other hand, very low frequencies are uninteresting because they cause the bubbles to evolve quasi-statically. 
Hence, without significant loss of generality, the dimensionless frequencies of $C_p$ are in the interval $[1/10,1/5]$. 
The phases of the waveforms that make up $C_p$ are independently sampled from a uniform distribution $\mathcal{U}$ in $[0,2\pi]$
\begin{gather} 
    \phi_i \sim \mathcal{U}([0,2\pi]), \quad i = 1, 2, ..., 6.
\end{gather}

Applications dictate the possible observed pressure amplitudes.
For example, significantly low pressures are not relevant for most applications. 
To set an empirical threshold approximating this condition, the pressures must not cause the used Monte Carlo simulation configuration to become numerically unstable.
The solver itself uses an adaptive 3rd-order Runge--Kutta scheme with minimum time-step $10^{-6}$ and relative error tolerance of $10^{-7}$. 
Thus, we design a pressure distribution from which all samples are numerically stable and physically realistic.
Algorithm~\ref{alg:forcing} details this process.

\begin{algorithm}[H]\label{alg:forcing}
\SetAlgoLined
$\alpha_i \sim \mathcal{U}([0,1]), \quad i = 1, 2, ..., 6$\;
$\alpha \leftarrow \sum_{i=1}^6 \alpha_i$\;
$\alpha_i \leftarrow \max \Big( 5 \alpha/ 3 \Big) \alpha_i, \quad i = 1, 2, \dots, 6$
 \caption{Forcing Amplitude Sampling}
\end{algorithm}

Previous experimental works can also be used to justify that the forcing constraint in algorithm~\ref{alg:forcing} avoids abrupt cavitation.
This is estimated by the cavitation number
\begin{gather}\label{Cavitation_number}
\begin{split}
    \sigma = \frac{1 - p_V(T_{0}) / p_{0}}{\rho U_{0}^2 / 2 p_{0}},
\end{split}
\end{gather}
where $p_V(T_{0})$ is the vapor pressure of the liquid at reference temperature $T_{0}$ and $U_0$ is the reference velocity~\cite{brennen2014cavitation}.
If the liquid cannot withstand negative pressures then vapor bubbles appear when the liquid pressure is $p_V$.
Thus, nucleation occurs when 
\begin{gather}
    \sigma \geq - \min_t \left\{ \frac{C_p(t)-1}{\rho U_0^2 / 2 p_0} \right\}.
\end{gather}
Without loss of generality, we can choose $\rho U_0^2/(2p_0) = 1$ to simplify the following calculations. For flows around axisymmetric headforms, with Reynolds number in the range of $4.5 \times 10^5$, if $\sigma \leq 0.40$, the formed nuclei grow explosively up to a certain bubble size~\cite{ceccio1991observations}. This phenomenon renders numerical simulations for flows with $\sigma \leq 0.40$ considerably more expensive compared to $\sigma > 0.40$ cases in order to achieve the same numerical accuracy. For the pressure profiles presented here, the case $\sigma \leq 0.40$ is avoided when using algorithm~\ref{alg:forcing}.

\begin{figure}[H]
    \centering
    \includegraphics[width=1.00\textwidth]{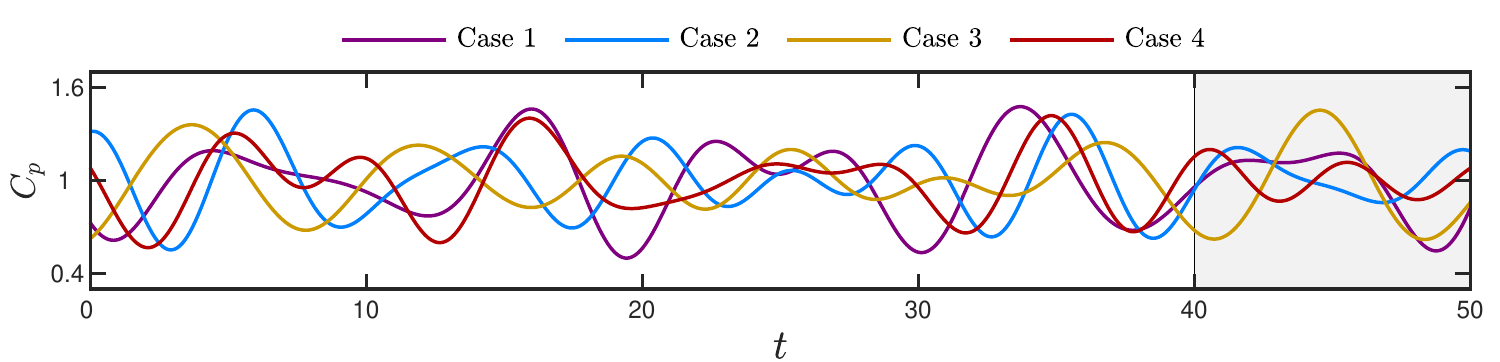}
    \caption{
        The time-history of example realizations of $C_p$. 
        Comparisons of the time-history of the evolved moments $\bmom$ and target moments $\bmom_{\mathrm{Targ.}}$ between different numerical schemes are performed in the shaded time-interval $t \in [40,50]$.
    }
    \label{fig:Pressure_Profiles}
\end{figure}

Figure~\ref{fig:Pressure_Profiles} shows example pressure profiles $C_p(t)$ that are used to test the fidelity of the hybrid moment inversion method.
Herein, the end of this time window, $t \in [40,50]$, is used to assess model fidelity. 
This enables the bubble dynamics to evolve from a specific initial state to one more representative of those found in actual flows.

\subsection{LSTM RNN training procedure}

We simulate $1000$ samples of individual bubbles for each realization of $C_p$.
The bubbles are initialized via samples from normal distributions with variances $\sigma_R^2$ and $\sigma_{\Rdot}^2$ for $R$ and $\Rdot$ respectively. 
Each case is evolved until $t = 50$, which in this dimensionless system corresponds to $50$ natural periods of bubble oscillations.
The individual bubble dynamics are then averaged to obtain the Monte Carlo reference statistics for each $C_p$ realization. 

For the numerical investigation, $200$ samples of $C_p$ from~\eqref{CP_Forcing} are used.
From those, $50$ are randomly selected for training, with the remaining $150$ cases used during testing. The number of samples used for training is chosen so that it is large enough to avoid over-fitting but small enough to still allow for the sampling of qualitatively different pressure profiles during testing.

\begin{figure}[H]
    \centering
    \subfloat[4-node]{{\includegraphics[width=0.45\textwidth]{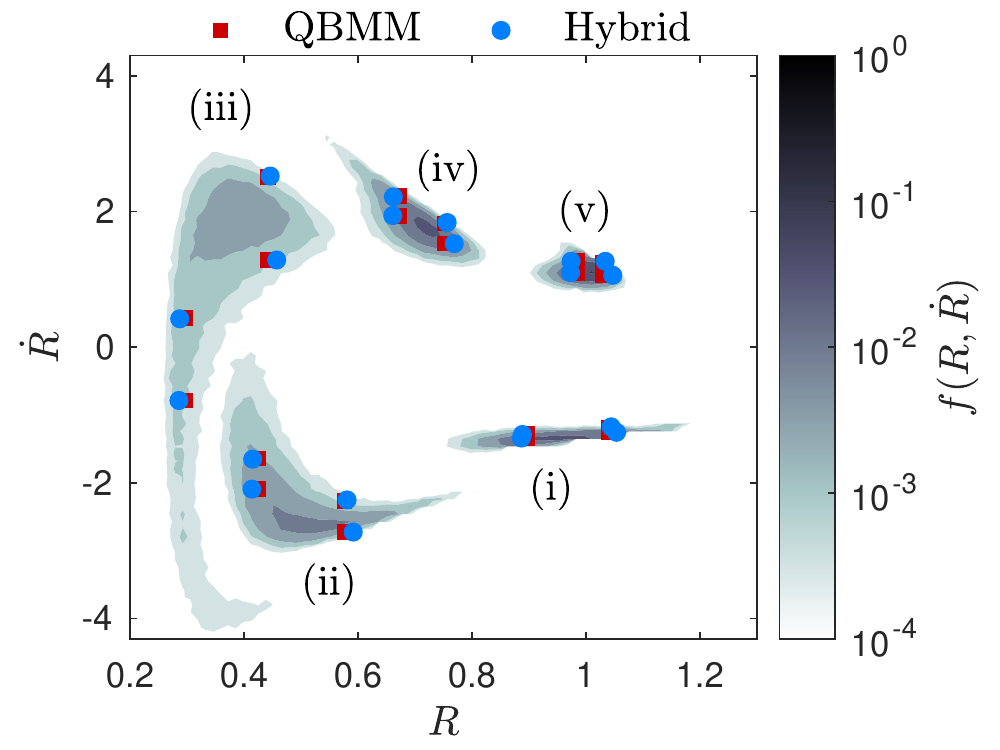} }}
    \subfloat[5-node]{{\includegraphics[width=0.45\textwidth]{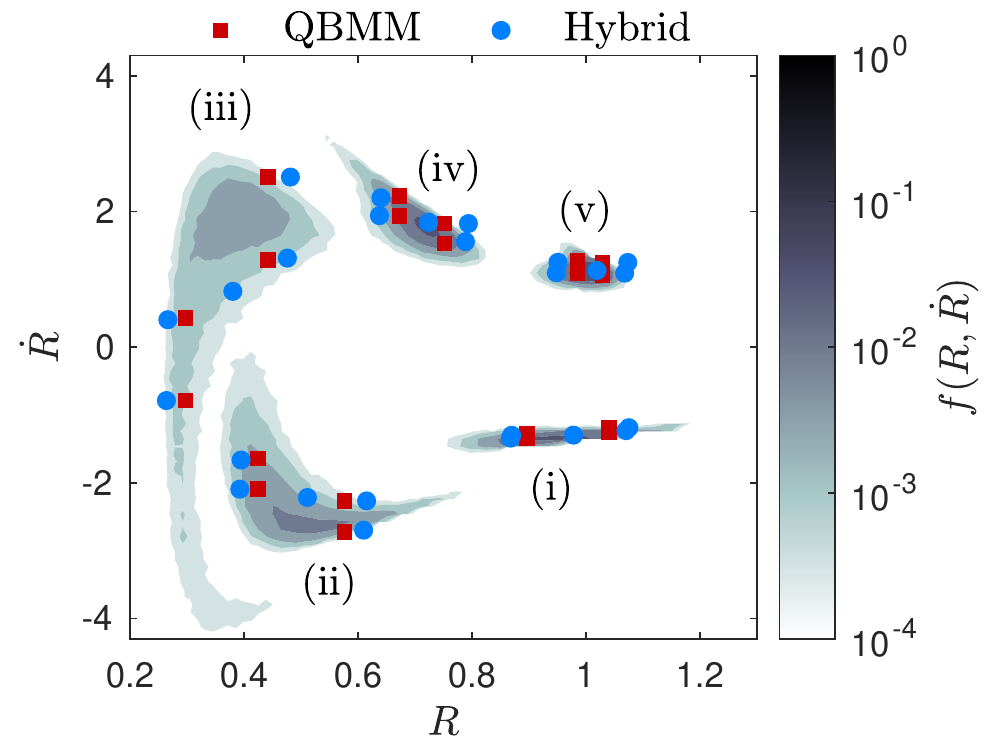} }}
    \caption{
        Temporal snapshots of $f$ computed via Monte Carlo and the positions of the quadrature nodes for the $4$-node CHyQMOM quadrature scheme (QBMM) and the $4$- and $5$-node hybrid CHyQMOM quadrature schemes (Hybrid).
        The labels (i)--(v) correspond to times $t = 43.9, 44.1, 44.2, 44.4$, and $44.6$, respectively.}
    \label{fig:Weight_Plots}
\end{figure}

Figure~\ref{fig:Weight_Plots} shows $f$ and the $\bxi$ for one realization of $C_p$ at different time instances. 
The same figure displays the CHyQMOM nodes as estimated by both the standard $4$-node CHyQMOM scheme and the $4$- and $5$-node hybrid schemes. 

During traning, the LSTM memory size is set to $256$ time-instances, with each $\delta t = 0.01$ apart. 
The Adam method~\cite{kingma2014adam} trained each neural network for $500$ epochs, minimizing the loss function~\eqref{Loss_Function}. A table with the values of the hyperparameters of the neural networks is presented in appendix~\ref{appendix:hyperparams}. A 4th-order Runge--Kutta adaptive time stepper evolves the predictions of the hybrid scheme.
The time integration scheme is adaptive, but the neural network predictions are uniform, so the neural network corrections are limited to the associated fixed time step $\delta t = 0.01$.

\subsection{Low-order moment evolution and error quantification}

\begin{figure}[H]
    \centering
    \includegraphics[width=1.00\textwidth]{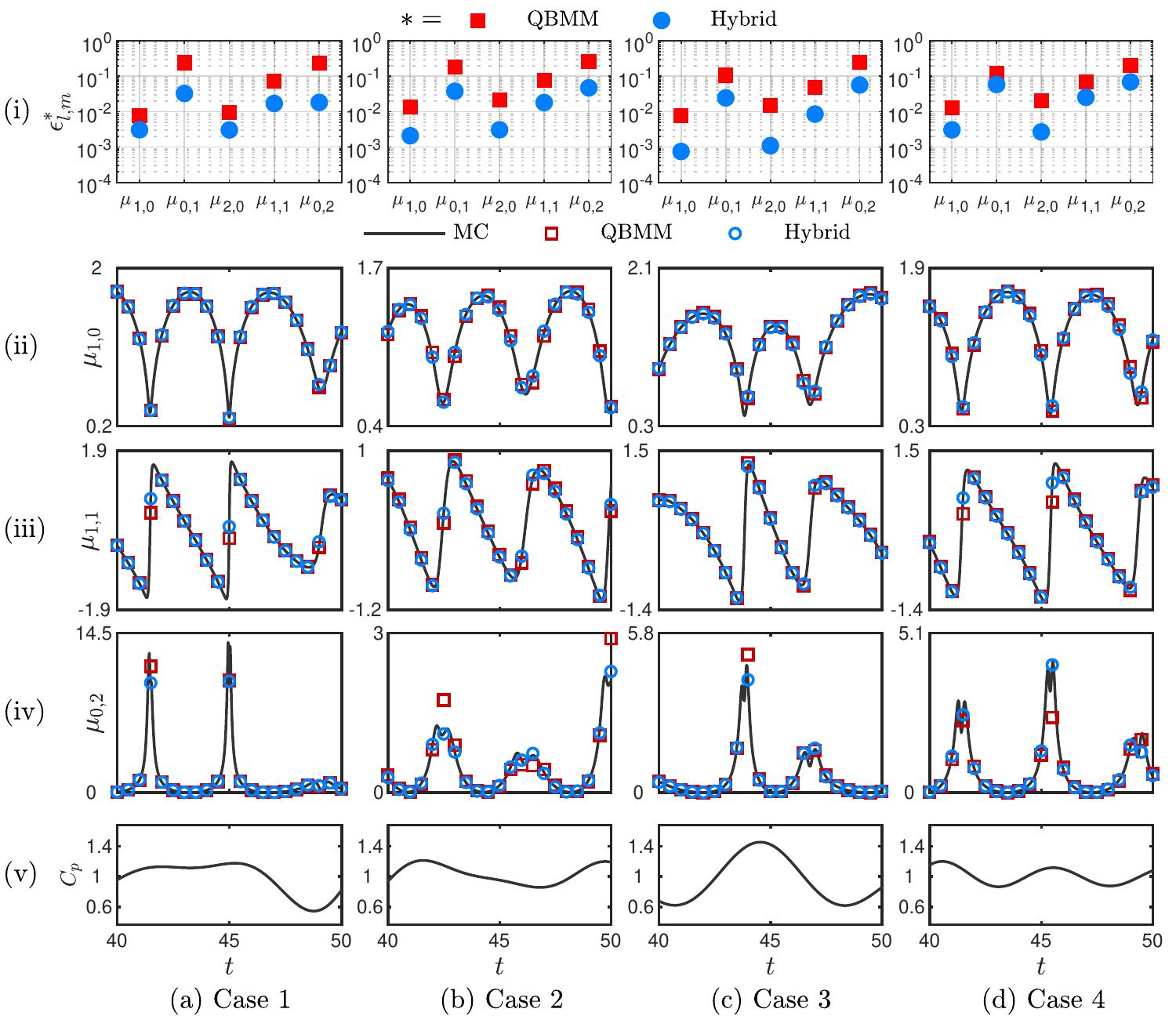}
    \caption{
        Low-order moment evolution for $4$-node CHyQMOM and hybrid CHyQMOM methods.
        Results compare against surrogate-truth Monte Carlo (MC) data.
    }
    \label{fig:LM_Results}
\end{figure}

The model-form relative error is computed via a discrete $L_2$ error
\begin{gather}
    \epsilon^{(\ast)}_{l,m} = 
    \sqrt{ \frac{\sum_{i=1}^{N_t} \left[ \mom_{l,m}^{\mathrm{(\ast)}}(t_i) - \mom_{l,m}^{\mathrm{(MC)}}(t_i) \right]^2 }{ \sum_{i=1}^{N_t} \left[\mom_{l,m}^{\mathrm{(MC)}}(t_i) \right]^2} },
\end{gather}
where (MC) indicates Monte Carlo surrogate truth data and $\ast = \{ \text{ML}, \text{QBMM} \}$.
The $t_i$ are $N_t = 5000$ uniformly spaced times in the interval $t \in [0,50]$. Results regarding the low-order moments are presented in figure~\ref{fig:LM_Results}.

Figure~\ref{fig:LM_Results}~(i) shows $\epsilon_{l,m}^\ast$ for the first- and second-order moments $\bmom$ for the $4$-node schemes.
Rows~(ii--iv) shows the evolution of select moments for $t \in [40,50]$ and row~(v) shows the corresponding $C_p$.
All results correspond to $4$ randomly selected testing samples (a)--(d) as labeled. We observe a smaller $\epsilon$ for the hybrid scheme than standard CHyQMOM for all $4$ cases considered. 
The largest errors for both approaches appear for moment $\mom_{0,2}$ (row~(iv)), which exhibits large and intermittent fluctuations when the bubbles collapse.  
The CHyQMOM method deviates most from the MC surrogate-truth data during intervals of high compression (small $C_p$), with hybrid CHyQMOM tracking the reference solution more accurately. 
Thus, the observed increase in accuracy varies significantly from case-to-case and moment-to-moment, from $10$ times smaller errors to only $20\%$ improvements. 

\subsection{High-order moment extrapolation}

\begin{figure}[H]
    \centering
    \includegraphics[width=1.00\textwidth]{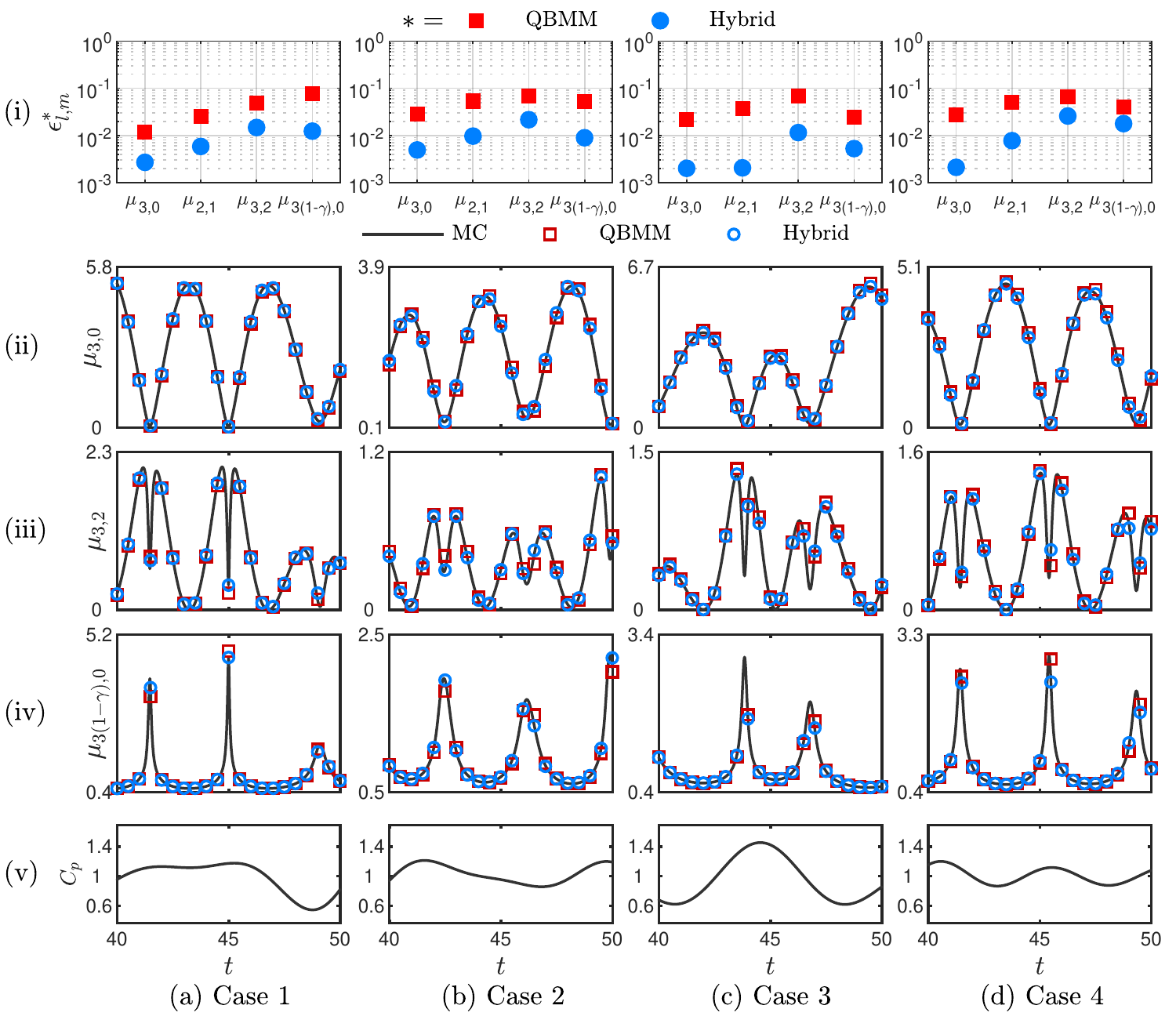}
    \caption{
        Comparison of predictions of target moments $\bmom_{\mathrm{Targ.}}$ for Monte Carlo data (black line), CHyQMOM predictions (red dashed line), and hybrid CHyQMOM predictions (blue line) for $4$ quadrature nodes.
    }
    \label{fig:HM_Results}
\end{figure}

The next quantity of interest is the $L^2$-error in predicting the target higher-order moments $\bmom_{\mathrm{Targ.}}$.
Figure~\ref{fig:HM_Results} presents these results for the same $4$ testing samples presented in figure~\ref{fig:LM_Results}. 
For all moments (ii)--(iv), hybrid CHyQMOM significantly improves the predictions of $\bmom_{\mathrm{Targ.}}$. 
This improvement is associated with the more accurate evolution of the low-order moments $\bmom$ and these targeted moments in the neural network training procedure.

To  better study the typical reduction in error for the $4$-node hybrid CHyQMOM scheme, we compute the percent improvement of the $L_2$ error as
\begin{gather}
    \mathcal{Q}_{l,m} = 100 \frac{\epsilon_{l,m}^{(\text{QBMM})} -\epsilon_{l,m}^{(\text{ML})}  }{ \epsilon_{l,m}^{(\text{QBMM})} }.
\end{gather}

\begin{figure}[H]
    \centering
    \includegraphics[width=1.00\textwidth]{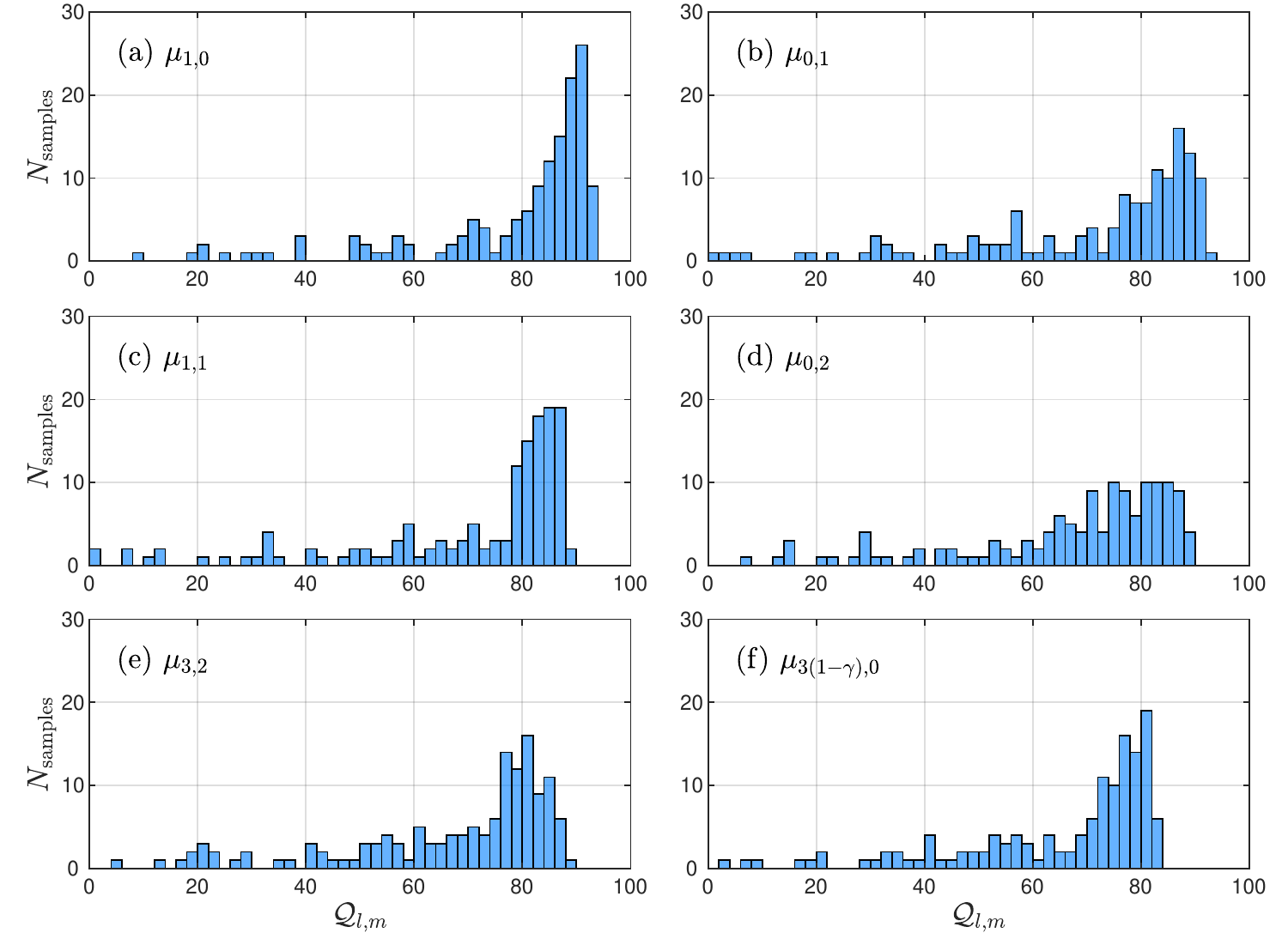}
    \caption{
        Histogram of $L^2$-error improvement $\mathcal{Q}$ for hybrid CHyQMOM over traditional CHyQMOM for example (a--d) low-order moments and (e--f) target high-order moments.
        Cases are drawn from 150 realizations of $C_p$. 
    }
    \label{fig:Histogram_Weights4}
\end{figure}

Figure~\ref{fig:Histogram_Weights4} shows a histogram of $\mathcal{Q}$  calculated for example low-order moments $\bmom$ and target moments $\bmom_{\mathrm{Targ.}}$.
For all testing samples, the $4$-node hybrid scheme improves the accuracy of the standard CHyQMOM method.

Further, for both low-order moments ($\bmom$) and target moments ($\bmom_{\mathrm{Targ.}}$), the error is reduced by more than $50\%$ for more than $80\%$ of the sampled $C_p$ cases.
The variation in error improvement is due to the amplitude range of the sampled $C_p$ and how closely the time-series of $C_p$ corresponds to one of the training samples.
Thus, results can improve by including more training samples.

\subsection{Additional quadrature nodes}

Another potential route to method improvement is to increase the number of quadrature nodes.
While the number of quadrature nodes can change, the evolved moments remain $\bmom$.
The algorithm for this is included in appendix~\ref{appendix:inversion}.

To quantify the effect of this change, $\epsilon^{1/2}_{l,m}(\ast)$ is computed, which is the median $\epsilon^{(\ast)}_{l,m}$ error among the $150$ test samples. 
We then define
\begin{gather}
    C_{l,m} = \frac{\epsilon^{1/2}_{l,m}(\text{ML}) }{\epsilon^{1/2}_{l,m}(\text{QBMM}) }
\end{gather}
to quantify the decrease in the $L^2$-error when using higher-node-count hybrid CHyQMOM compared to the standard $4$-node CHyQMOM.

\begin{figure}[H]
    \centering
    \includegraphics[width=1.00\textwidth]{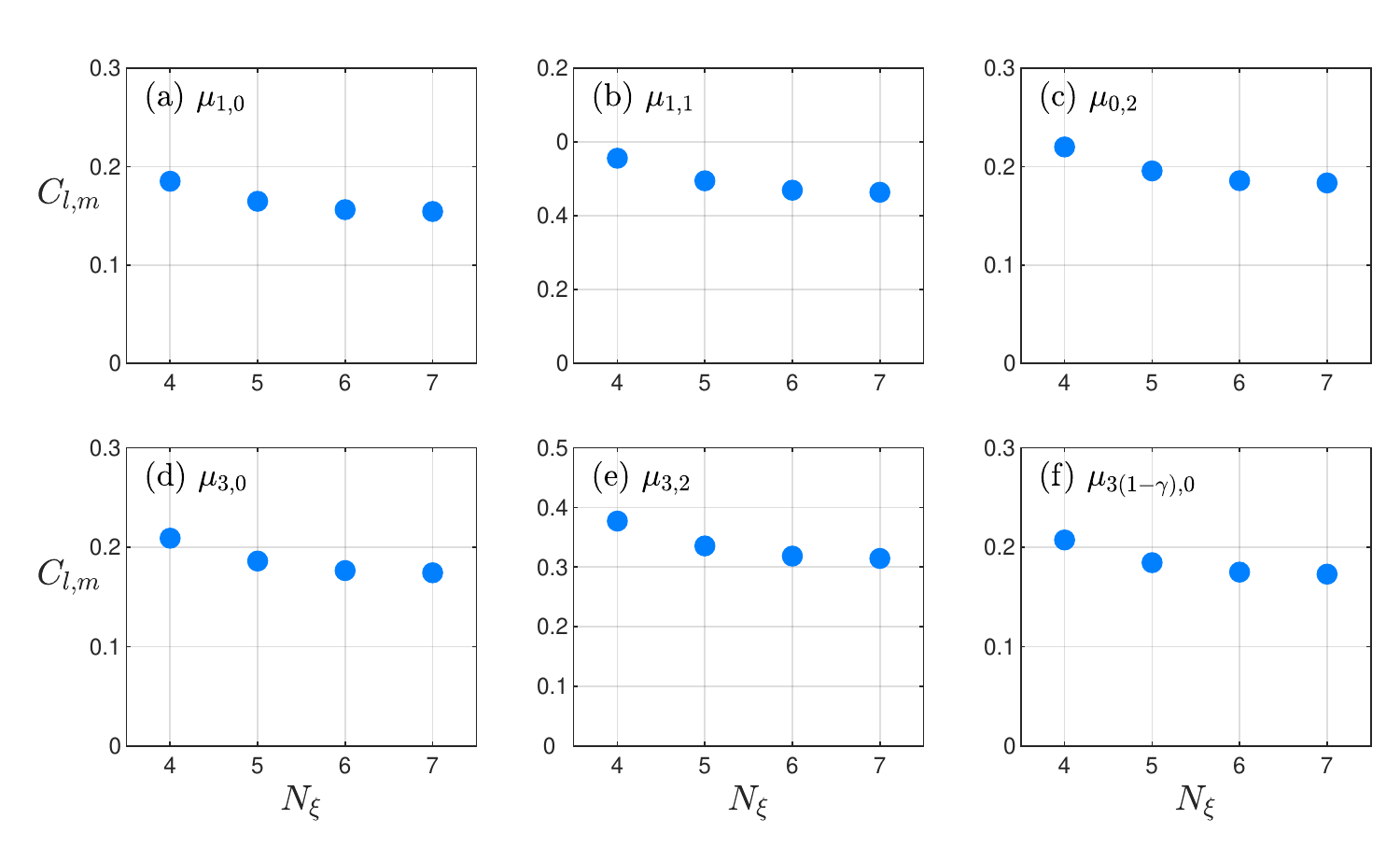}
    \caption{Median error decrease while using hybrid CHyQMOM over $4$-node CHyQMOM for different numbers of nodes for the hybrid CHyQMOM scheme.}
    \label{fig:Median_Weights4}
\end{figure}

Figure~\ref{fig:Median_Weights4} shows that the accuracy of the hybrid predictions is improved as the number of nodes $N_\xi$ increases.
However, the error improvement gains diminish once $7$ nodes are reached. 
Further, including additional nodes to the quadrature rule increases the computational time needed to perform a single time-step evolution for the system. 
The computational cost of $4$-node hybrid CHyQMOM per time-step is $8.9$ times the cost of CHyQMOM.
For $5$, $6$, and $7$ nodes, the hybrid method costs per time-step are $11.5$, $13.8$, and $16.2$ times that of $4$-node CHyQMOM\footnote{
These simulations were performed using $\mathtt{PyQBMMlib}$~\cite{bryngelson2020qbmmlib} on a single core of a $\SI{2.3}{\giga\hertz}$ Intel Core i$9$ CPU.}.
Hence, diminishing improvements are observed as the number of nodes increases to more than $6$, as the simulations require significantly more computations per time-step for comparable accuracy. 

\section{Conclusions}\label{s:conclusions}

A data-informed conditional hyperbolic quadrature method for statistical moments was presented. 
The method was applied to the statistics of a population of spherical bubbles oscillating in response to time-varying randomized forcing. 
The forcing is designed to resemble any possible function with a banded frequency spectrum from $1/5$ to $1/10$ the natural frequency of the bubbles. 
Results showed that the method reduces closure errors when compared against a standard $4$-node CHyQMOM scheme. 
The hybrid method reduced errors more significantly for the extrapolated higher-order moments that close the phase-averaged flow equations.
This improvement was achieved via recurrent neural networks that include time history during training.
This result is significant because higher-order QBMM schemes are generally numerically unstable for this problem, so another approach is required to improve accuracy.
Thus, while the presented hybrid scheme is about a factor of $10$ more expensive than CHyQMOM, its numerical cost should be viable for many applications. \vskip6pt

\enlargethispage{20pt}



\aucontribute{
All authors jointly conceived the problem formulation.  A.C and S.B. designed the study.  A.C.\ developed the ML code, conducted the numerical experiments, and drafted the manuscript. S.B.\ contributed the bubble dynamics code and further developments of $\mathtt{PyQBMMlib}$.  All authors reviewed and edited the manuscript.
}

\competing{We declare we have no competing interests.}

\funding{
The US Office of Naval Research supported this work under grant numbers N0014-17-1-2676 and N0014-18-1-2625.
Computations were performed via the Extreme Science and Engineering Discovery Environment (XSEDE) under allocation CTS120005, supported by National Science Foundation grant number ACI-1548562.
}

\ack{The authors thank Rodney Fox and Alberto Passalacqua for valuable discussion of quadrature moment methods.}




\bibliography{main_article_RSTA} 
\bibliographystyle{abbrv}


\appendix
\section{CHyQMOM inversion algorithm}\label{appendix:inversion}

This is the inversion algorithm for the $4$-node CHyQMOM scheme. Given the first and second-order moments $\{\mom_{1,0},\mom_{0,1},\mom_{2,0},\mom_{1,1},\mom_{0,2}\}$ it computes the nodes $(\xi_{i},\dot{\xi}_{i})$ and corresponding weights $w_i$ for $i = 1, 2, 3, 4$, in phase-space. In this work we assume $\mom_{0,0} = 1$. To tail the algorithm to our hybrid scheme of arbitrary number of quadrature nodes, the algorithm adds some fictitious extra nodes to the scheme with zero-valued weights to match the desired number of nodes of the hybrid scheme.

\begin{algorithm}[H]
\SetAlgoLined
$w_i = 0.25, \quad 1 \leq i \leq 4$.\;
$w_i = 0.00, \quad 4 < i \leq N$\;
$\sigma_R = \sqrt{\mom_{2,0} -\mom_{1,0}^2}$\;
$\alpha = \frac{ \mom_{1,1} -\mom_{1,0} \mom_{0,1} }{ \sigma_R }$\;
$\sigma_{\Rdot} = \sqrt{ \mom_{0,2} -\alpha^2 -\mom_{0,1}^2 }$\;
$\xi_i =\mom_{1,0} +\sigma_R$\; 
$\xi_2 =\mom_{1,0} +\sigma_R$\; 
$\xi_3 =\mom_{1,0} -\sigma_R$\; 
$\xi_4 =\mom_{1,0} -\sigma_R$\; 
$\xi_i =\mom_{1,0}, \quad 4 < i \leq N$\; 
$\dot{\xi}_1 = \mom_{0,1} +\alpha +\sigma_{\Rdot}$\;
$\dot{\xi}_2 = \mom_{0,1} +\alpha -\sigma_{\Rdot}$\;
$\dot{\xi}_3 = \mom_{0,1} -\alpha +\sigma_{\Rdot}$\;
$\dot{\xi}_4 = \mom_{0,1} -\alpha -\sigma_{\Rdot}$\;
$\dot{\xi}_i = \mom_{0,1}, \quad 4 < i \leq N$\;
\caption{CHyQMOM inversion algorithm}
\end{algorithm}

\section{Neural network hyperparameters}\label{appendix:hyperparams}

\begin{table}[H]
\centering
\begin{tabular}{  r | l   }
 \hline\hline
 \textbf{Hyperparameter} & \textbf{Value}\\
 \hline
 Epochs   & $500$\\
 Learning rate & $10^{-5}$\\
 Batch size & $32$\\
 Activation function & tanh\\
 Recurrent activation function & hard sigmoid\\
 Dropout coefficient & $0.10$\\
 Recurrent dropout coeffificent & $0.10$\\
 LSTM is stateful & True\\
 Kernel initializer & Zeros\\
 Recurrent initializer & Zeros\\
 Bias initializer & Zeros\\
 \hline\hline
\end{tabular}
\caption{Hyperparameters used to train the neural networks.}
\label{t:hyper}
\end{table}

\end{document}